\documentclass[12pt]{article}
\usepackage{graphicx}
\usepackage[all,poly]{xy}
\usepackage{url}
\title{
Alternative parametrizations and reference priors for decomposable
discrete graphical models}
\author{
Guido Consonni \\
Dipartimento di Economia Politica e Metodi Quantitativi\\  University of Pavia\\ Via San Felice 5,
27100 Pavia, Italy\\
\url{guido.consonni@unipv.it}
\and   H\'el\`ene Massam\\
Department of
Mathematics and Statistics, York University \\ Toronto, M3J 1P3,
Canada\\
\url{massamh@mathstat.yorku.ca}
}
\date{ }
\usepackage{latexsym}
\begin{document}
\maketitle

\def\reel{\hbox{I\hskip-2pt R}}
\def\comp{\hbox{I\hskip-6pt C}}
\def\quat{\hbox{I\hskip-2pt H}}
\def\gene{\hbox{I\hskip-2pt D}}
\def\natu{\hbox{I\hskip-2pt N}}
\def\tr{\;\mathrm{tr}\;}
\def\<{\langle}
\def\>{\rangle}
\def\P{{\bf P}}
\def\E{{\bf E}}
\def\n{{\boldmath n}}
\def\b{{\boldmath b}}
\def \m{{\bf m}}
\def\z{{\boldmath z}}

\def\<{\langle}
\def\>{\rangle}

\def\th{th\'eor\`eme\ }
\def\ths{th\'eor\`emes\ }
\def\reel{\hbox{I\hskip -2pt R}}
\def\pe{\hbox{I\hskip -2pt P}}
\def\u{{\bf u}}
\def\v{{\bf v}}
\def\esp{{\bf E}}
\def\<{\langle}
\def\>{\rangle}
\def\x{{\bf x}}
\def\y{{\bf y}}
\def\f{{\bf f}}
\def\p{{\bf p}}
\def\w{{\bf w}}
\def\s{{\bf s}}
\def\W{{\bf W}}
\def\Z{{\bf Z}}
\def\Q{{\bf Q}}
\def\tr{\textmd{trace}\,}
\def\grad{\textmd{grad}\,}

\def\ix{\smallskip \hangindent=1.5em \hangafter=1 \noindent}

\textheight =22.5 cm
\textwidth =15 cm
\voffset =-0.5 in
\hoffset =0 in
\headheight =0 cm
\renewcommand{\thesection}{\arabic{section}}
\renewcommand{\theequation}{\arabic{section}\mathrm{.}\arabic{equation}}

\newtheorem{theorem}{Theorem}[section]
\newtheorem{definition}{Definition}[section]
\newtheorem{subsecdef}{Definition}[subsection]
\newtheorem{lemma}{Lemma}[section]
\newtheorem{prop}{Proposition}[section]
\newtheorem{remarks}{Remarks}[section]
\newtheorem{remark}{Remark}[section]
\newtheorem{example}{Example}[section]
\newtheorem{perty} {Property}[section]
\newtheorem{cor} {Corollary}[section]

\newenvironment{pff}{\hspace*{-\parindent}{\bf Proof:}}{\hfill $\Box$
\vspace*{0.2cm}}
\newenvironment{pff.}{\hspace*{-\parindent}{\bf Proof}}{\hfill $\Box$
\vspace*{0.2cm}}
\renewcommand\theequation{\thesection .\arabic{equation}}
\let\subs\subsection

\let\sect\section
\renewcommand\section[1]{\setcounter{equation}{0}\sect{#1}}

\newcommand{\pperp}{{\ \perp\!\!\!\perp\ }}
\def\tr{\;\mathrm{tr}\;}
\def\Si{\Sigma^{\cal V}}

\newcommand{\cK}{{\cal K}}
\newcommand{\cM}{{\cal M}}
\newcommand{\IR}{{\Bbb R}}
\newcommand{\IC}{{\Bbb C}}
\newcommand{\ID}{{\Bbb D}}
\newcommand{\IO}{{\Bbb O}}
\newcommand{\IN}{{\Bbb N}}
\newcommand{\IH}{{\Bbb H}}
\vspace{-0.2cm}
\begin{abstract}

For a given discrete decomposable graphical model, we identify
several alternative parametrizations, and construct the
corresponding reference priors for suitable groupings of the
parameters. Specifically, assuming  that the cliques of the graph
are arranged in a perfect order, the parameters we consider
are conditional probabilities of clique-residuals given
separators, as well as generalized log-odds-ratios.
We also consider
a parametrization associated to a collection of variables representing
a cut for the statistical model. The reference priors we obtain do
not depend on the order of  the groupings, belong  to a conjugate
family, and are proper.

\vspace{.5cm}

\textit{Some key words}:  Clique; Conjugate family;  Contingency table; Cut;
Log-linear model;  Multinomial model; Natural exponential family.
\vspace{.5cm}

\end{abstract}

\section{Introduction}
\label{sec:intro}
Graphical models, see e.g. Lauritzen (1996), are statistical models such that  dependencies between variables are
expressed by means of a graph. The study of graphical models is an established and active area of applied
and theoretical research.
Directed graphs for  discrete variables, often called  Bayesian
networks, see e.g. Cowell et al. (1999), have been used in a variety of applied domains, and represent
the engine of probabilistic expert systems.
On the other hand, undirected graphical models for the  analysis of  discrete data
   are best employed for the analysis of  multi-way contingency tables, and
represent  a useful  subset of hierarchical log-linear models .

In this paper we are concerned with  the Bayesian analysis of discrete undirected graphical
models, whose underlying graph is decomposable.
When working in a Bayesian framework, a prior distribution on the parameter space is required.
Priors for undirected discrete graphical, or more generally, log-linear models
have been considered in Dawid and Lauritzen (1993), Madigan and York (1995), Dellaportas and
Forster (1999),  Kings and Brooks (2001), Dellaportas and Tarantola (2005).

Despite the adoption of  reasonably simplified models,
 prior elicitation  still represents a major
 concern even for moderately large graphs, because of the very
 high number of parameters involved.
This naturally suggests to search for  default, or objective, priors,
requiring a minimal subjective input and essentially model-based.
 However there is now evidence, see e.g. Berger (2000) and Casella (1996),
that naive approaches
 based on flat non-informative priors are largely inadequate in multi-parameter settings.
In this context, reference analysis  provides one of the most successful general methods to
derive default  prior distributions. For a recent and informative review see
Bernardo (2005). While  the algorithmic complexity for the construction of reference priors can be
substantial, it is known that suitable re-parametrizations of the model may considerably simplify
the task, see for instance
Consonni et al. (2004) and Consonni and Leucari (2006).

We  address two specific issues in this paper:
identifying alternative parametrizations for a given discrete graphical model,
 and constructing the
corresponding reference priors.
More precisely, in \S \ref{sec:parametrizations}
we consider several  parametrizations: conditional
probabilities of clique-residuals given separators, as well as
generalized log-odds ratios  that arise as canonical
parameters of equivalent exponential family representations of the underlying sampling
distribution, and explicate their mutual relationships. In \S \ref{sec:reference priors} we provide
the expressions for the corresponding reference priors, and discuss  their main properties.
In \S 4 we present a parametrization
associated to a cut in the graphical model and derive the corresponding reference prior.
Some points for discussion are summarized in the last section. Technical details
for the proof of the relationships between various parametrizations are given in the Appendix.


\section{Generalized log-odds-ratios parametrizations}
\label{sec:parametrizations}
\subsection{Preliminaries}
Let us recall some basic facts about undirected graphs and
graphical models: for further details the reader is referred to
Lauritzen (1996, ch. 2). An undirected graph $G$ is a pair $(V,E)$
where $V$ is a finite set of vertices and $E$ the set of
edges, an edge being an unordered pair $\{ \gamma,\delta\}, \, \gamma \in V,\delta \in
V, \gamma \neq \delta$. Henceforth the graph $G$ is assumed to be decomposable.
For a given ordering
 $C_1,\ldots,C_k$  of the cliques, we will use the following notation
 $$H_l=\cup_{j=1}^l C_j, l=1,\ldots,k,\;\;S_l=H_{l-1}\cap C_l,\; l=2,\ldots,k,\;\;R_l=C_l\setminus S_{l},\; l=2,\ldots,k\;.$$
 A given ordering of the cliques is said to be perfect if for any $l >1$ there is an $i < l$ such that $S_l \subseteq C_i$. When we have a perfect ordering of the cliques, the $S_l, l=2,\ldots,k$ are minimal separators. The $H_l$ and $R_l$ are called respectively, the $l$-th history and $l$-th residual.

A graphical model, Markov with respect to a given graph $G$, is a family of
probability distributions on $(X_{\gamma}, \, \gamma \in V)$ such that
$X_{\delta}$ is independent of $X_{\gamma}$ given $X_{V\setminus \{\delta,\gamma\}}$ whenever
 $\{ \gamma, \delta \}$ is not in $E$.

In this paper we shall focus on contingency tables
arising  from the classification of $N$ units according to a finite
set $V$ of criteria, see Lauritzen (1996, Ch. 4).
  Each criterion  is represented by
a variable $X_{\gamma}$, $\gamma\in V$, which takes values in a finite set ${\cal
I}_{\gamma}$. Let ${\cal I}=\times_{\gamma\in V}{\cal I}_{\gamma}$.
The cells of the table are the elements
\begin{equation}
i=(i_{\gamma},\;\gamma\in V), \;i\in {\cal I}.
\end{equation}
Each of $N$ individuals falls into cell $i$
independently with a probability $p(i)$; we let $p=(p(i), \; i \in {\cal I})$, with
$\sum_{i \in {\cal I}}p(i)=1$.
 Furthermore, we write $n(i)$ for the $i$-th cell-count and $n=(n(i), \; i \in {\cal I})$, with
 $\sum_{i \in {\cal I}}n(i)=N$.

 We consider here the model
 ${\cal M}_G$, which, for a given $G$ and a given integer $N$,
 is the set of multinomial ${\cal M}(N,p)$ distributions with $N=\sum_{i\in{\cal I}}n(i)$ and
$p=(p(i), \; i \in {\cal I})$ in the $|{\cal I}|-1$ dimensional simplex, which are Markov with respect to $G$.
\vspace{2mm}

From now on, we adopt the notation  \lq
\lq $D\subseteq_0 V$\rq \rq{} to mean that $D$ may be the empty
set while \lq
\lq $D\subseteq V$\rq \rq{} excludes the empty
set.
Let ${\cal E}$ denote the power set of $V$, excluding the empty set, i.e.
\begin{eqnarray*}
{\cal E}=\{F \subseteq V, F\not =\emptyset \}.
\end{eqnarray*}

 For $D\in {\cal E}$,
\begin{equation}
\label{basei} i_D=(i_{\gamma},\;\gamma\in D), \;\;\mbox{and}\;\;
n(i_D),\;\;i_D\in {\cal I}_D=\times_{\gamma\in D}{\cal I}_{\gamma}
\end{equation}
 denotes a cell in the $D$-marginal table, and its corresponding count.
 We therefore have
 \begin{equation}
 \label{niD}
 n(i_D)=\sum_{j\in {\cal I}|j_D=i_D}n(j)=\sum_{j_{V\setminus D}\in {\cal I}_{V\setminus D}}n(i_D,j_{V\setminus D}).
 \end{equation}
 Note that $n(i_{\emptyset})=N$.
 For $F,D$ in ${\cal E}$,
 we use the notation $p^D(i_D)$ and $p^{D|i_F}(i_D)$ to denote, respectively,
 the marginal and the conditional probabilities
 \begin{eqnarray}
 p^D(i_D)
 &=&\sum_{j\in {\cal I}|j_D=i_D}p(j)\\
 p^{D|i_F}(i_D)
&=&\frac{p^{D\cup F}(i_D,i_F)}{p^F(i_F)}.
 \end{eqnarray}
Assuming that \lq \lq $0$\rq \rq{} indicates one of the levels for each variable,
   we let $i_{\gamma}^*$ denote the\lq \lq $0$\rq \rq-level in ${\cal I}_{\gamma}$, so that
$$i^*=(i^*_{\gamma},\;\gamma\in V)$$
denotes the cell with all components equal
to $0$.

\begin{definition}
For $D\in {\cal E}$, we define
\begin{equation}
{\cal I}_{D}^*=\{i_D\;|\; i_{\gamma}\not =i^*_{\gamma},\forall
\gamma\in D\}.
\end{equation}
\end{definition}
In words, ${\cal I}_{D}^*$ is the set of marginal cells $i_D$ such that none of their components is equal to 0.
We set ${\cal I}_{V}^*={\cal I}^*$.
 For
example, if $D=\{a,b,c\}$, $a$ takes the values $\{0,1,2,3\}$, $b$
takes the values $\{0,1,2\}$, $c$ takes the values $\{0,1\}$, then
\begin{eqnarray*}
{\cal
I}_{D}^*&=&\{(1,1,1),(2,1,1),(3,1,1),(1,2,1),(2,2,1),(3,2,1)
\}
\end{eqnarray*}

%

\subsection{The saturated case}
We assume here that $G$ is complete and ${\cal M}_G$ is therefore the saturated multinomial
model for $n=(n(i),i \in {\cal I}))$
. The multinomial probability function  is usually written in terms of the cell probabilities $p=(p(i),i
\in {\cal I})$  as
\begin{eqnarray}
\label{saturated multinomial}
f(n|p)
= \frac{N!}{\prod_{i \in {\cal I}}n(i)!}\prod_{i\in {\cal
I}}p(i)^{n(i)},
\end{eqnarray}
where the only restriction
 on the parameters $p(i)$ is  $\sum_i p(i)=1$.
It is convenient  to regard the
multinomial coefficient in (\ref{saturated multinomial}) as being
part of the dominating measure, so that the actual density is simply
$\prod_{i\in {\cal I}}p(i)^{n(i)}$.
Assuming that all probabilities are positive,
 the density (\ref{saturated multinomial}), with respect to a suitable dominating measure,
  can be represented in exponential family form as
\begin{equation}
\label{multinn(i)}
\prod_{i\in {\cal I}}p(i)^{n(i)}=\exp \left \{ \sum_{i\in {\cal
I},i\not=i^*}n(i)\xi(i)-N\log\big(1+\sum_{i\in {\cal
I},i\not=i^*}e^{\xi(i)}\big) \right \},
\end{equation}
where
\begin{equation}
\xi(i)=\log \frac{p(i)}{p(i^*)},\;\;\;i\in {\cal I}^*
\end{equation}
are the usual log-odds, relative to the benchmark cell $i^*$.
We recognize in (\ref{multinn(i)}) a natural exponential family (henceforth abbreviated NEF),
with canonical parameters $\xi(i)$ and canonical statistics $n(i)$, $i \neq i^* $.
For a review of NEFs, see e.g. Kotz, Balakrishnan, and Johnson (2000, ch. 54).

In this paper,
we shall work with
NEF-representations alternative to (\ref{multinn(i)}), featuring  different canonical
statistics and  their corresponding canonical parametrizations,
the latter representing various generalized log-odds-ratios of  joint probabilities,
residual-conditional probabilities or clique-marginal probabilities.
For the saturated model, we need consider only the generalized log-odds-ratio of model probabilities defined
as follows.

 \begin{definition}
For all $D\subseteq_0 V$ and $i_D\in {\cal I}^*_D$ we define the
log-linear parameters
\begin{equation}
\label{thetaid} \theta(i_D)=\log \prod_{F\subseteq_0
D}p(i_F,i^*_{V\setminus F})^{(-1)^{|D\setminus F|}},
\end{equation}
\end{definition}
 Note that for $F=\emptyset,\; p(i_F,i^*_{V\setminus F})=p(i^*)$ and
$\theta (i_{\emptyset})=\theta(i^*)=\log p(i^*).$ The parameters $\theta(i^*)$ and $p(i^*)$ are not free but functions of the other $\theta$ or $p$ parameters.
 We also
emphasize  the fact that while  $\theta(i_D)$ is indexed by the marginal cell $i_D$, it
is a function of the joint probabilities $p(i)$ in the full table.

Making the change of variables
$$\Big(n(i), i\in {\cal I}\setminus \{i^*\}\Big)\mapsto \Big(n(i_D), D\subseteq V, i_D\in {\cal I}^*_D\Big),$$ it is relatively easy to show the following expression of the multinomial distribution.
\begin{prop}
\label{prop:basicmult} The NEF-representation of the saturated multinomial
model in terms of the log-linear parameters $\theta(i_D)$ is given
by
 \begin{eqnarray}
\hspace{-.5cm}
\prod_{i\in {\cal I}}p(i)^{n(i)}
=
\exp \left \{\sum_{D\subseteq V} \sum_{i_D\in {\cal
I}^*_D}n(i_D)\theta(i_D)-N\log \big(1+\sum_{D\subseteq
V}\sum_{i_D\in {\cal I}^*_D}
 \exp{\sum_{F\subseteq D}\theta(i_F)}\big) \right \}\label{basicmult}.
\end{eqnarray}
\end{prop}

We remark that the canonical parameter $\theta(i_D)$ in
(\ref{thetaid}), $D \subseteq V$, is defined only for $i_D\in {\cal
I}^*_D$, i.e. all those cell-configurations having no component
equal to $0$; alternatively the remaining components indexed by
$i_D\in {\cal I} \setminus{\cal I}^*_D$ may be regarded as being set
to zero in (\ref{basicmult}), and thus  satisfy the usual \lq \lq
corner constraint\rq \rq{}  used for instance in GLIM.
Furthermore, the canonical statistics $n(i_D)$ represent the
marginal counts for all cells $i_D, \; D\subseteq V$ and $i_D\in
{\cal I}^*_D$.

 \subsection{The case for $G$ decomposable}
 If the multinomial model is Markov with respect to a given decomposable, non complete, graph $G$,
 it is a simple consequence of the Hammersley-Clifford theorem
 (see Lauritzen, 1996, p. 36 and Liu and Massam, 2007) that the
 model is Markov with respect to $G$ if and only if for
$i_D\in {\cal I}^*_D,\; D\subseteq V$
 \begin{equation}
 \label{thetaidequalzero}
 \theta(i_D)=0\;\;\mbox{whenever}\;\;D\;\;\mbox{is not complete in G}.
\end{equation}
 The model (\ref{basicmult}) satisfying (\ref{thetaidequalzero}) as
the multinomial  model ${\cal M}_G$ Markov with respect to $G$. More briefly, we refer to it as the multinomial Markov model.

For any subset $A\subseteq V$ of the vertex set, define
 \begin{eqnarray}
 {\cal D}^A &=& \{D\subseteq A\;|D\;\;\mbox{is complete}\} \\
{\cal D}_0^A &=& \{D\subseteq_0 A\;|D\;\;\mbox{is complete}\}.
 \end{eqnarray}
 To simplify notation, we will write ${\cal D}$ for ${\cal D}^V$.
 We are going to present in this subsection three parametrizations for ${\cal M}_G$.

 \noindent  The first parametrization
 is in terms
 of the log-linear parameters defined in (\ref{thetaid}) with canonical parameter
 \begin{equation}\label{thetamodel}
 \theta^{mod}=\theta({\cal D})=(\theta(i_D),\;D\in {\cal D},\;i_D\in {\cal I}^*_D)\;
 \end{equation}
 and corresponding canonical statistic
 \begin{equation}
 \label{statmodel}
 n({\cal D})=\big(n(i_D),\;D\in {\cal D},\;i_D\in {\cal I}^*_D\big).
 \end{equation}
 It will also be convenient to use the notation
 \begin{equation}
 k(\theta({\cal D}^A))=\log \big(1+\sum_{D\subseteq
A}\sum_{i_D\in {\cal I}^*_D}
 \exp{\sum_{F\subseteq D}\theta(i_F)}\big)
 \end{equation}
 for the cumulant generating function, and the  notation
 \begin{equation}
 \<\theta({\cal D}^A),n({\cal D}^A)\>=\sum_{D \subseteq A} \sum_{i_D\in {\cal
I}^*_D}\theta(i_D)n(i_D).
 \end{equation}
 for the inner product.
 The NEF representation  in terms of $\theta({\cal D})$
 can then be immediately derived from (\ref{basicmult}) as follows.
  \begin{prop}
\label{prop:markovmult} Let $G$ be a decomposable graph. The NEF-representation
of the multinomial Markov model in terms of the parametrization $\theta^{mod}$ is given by
 \begin{eqnarray}
 \label{markovmult}
\exp \{\<\theta({\cal D}),n({\cal D})\>-N\; k(\theta({\cal
D}))\}.
\end{eqnarray}
\end{prop}
\vspace{2mm}

\noindent Let us now introduce a second parametrization which is relative to the marginal distribution for $C_1$
and the conditional distributions for $R_l$ given $S_l$.
 For a given  perfect ordering $C_1,\ldots,C_k$ of the cliques of
 $G$, the Markov property implies (see Lauritzen, 1996, p. 90)
 \begin{equation}
p(i)=\frac{\prod_{l=1}^kp^{C_l}(i_{C_l})}{\prod_{l=2}^kp^{S_l}(i_{S_l})}.
 \end{equation} As a consequence we can write the multinomial
 density (\ref{saturated multinomial}) as
%
%
 \begin{eqnarray}
\prod_{i\in {\cal I}}p(i)^{n(i)}&=&\prod_{i\in {\cal
I}}\left(\frac{\prod_{l=1}^kp^{C_l}(i_{C_l})}{\prod_{l=2}^kp^{S_l}(i_{S_l})}\right)^{n(i)}=\prod_{i\in
{\cal I}}\left(p^{C_1}(i_{C_1})\prod_{l=2}^k
p^{R_l|i_{S_l}}(i_{R_l})\right)^{n(i)}\nonumber\\
&=&\prod_{i_{C_1}\in {\cal I}_{C_1}} (p^{C_1}(i_{C_1}))^{n(i_{C_1})}
\prod_{l=2}^k\prod_{i_{C_l}\in {\cal I}_{C_l}}\left(p^{R_l|i_{S_l}}(i_{R_l})\right)^{n(i_{C_l})}\nonumber\\
&=&\prod_{i_{C_1}\in {\cal I}_{C_1}} (p^{C_1}(i_{C_1}))^{n(i_{C_1})}
\prod_{l=2}^k\prod_{i_{S_l}\in {\cal I}_{S_l}}\prod_{i_{R_l}\in
{\cal
I}_{R_l}}\left(p^{R_l|i_{S_l}}(i_{R_l})\right)^{n(i_{C_l})}\;.\label{cliquesepdecomp}
\end{eqnarray}
Note that (\ref{cliquesepdecomp}) expresses the multinomial Markov model in terms of the marginal probabilities
in the
$C_1$-table,  as well as the conditional probabilities in the
$i_{S_l}$-slice of the $R_l$-table, for $l=2,\ldots,k$.

Formally, for $B \subset V$ and $A \subset V$ with $ A \cap
B=\emptyset$, the $i_B$-slice of the $A$-table is obtained by
classifying, according to the factors in $A$, only those units that
belong to the marginal $i_B$-cell (for the notion of \lq \lq
slice\rq \rq{} in a contingency table see Lauritzen, 1996, p. 68).

 Let us now define the log-linear parameters corresponding to the factorization (\ref{cliquesepdecomp}).
 \begin{definition}
 \label{def:thetamarg}
For each clique $C_l$, $l=1,\ldots,k$,
we define
\begin{equation}
\label{thetamarg} \theta^{C_l}(i_D)=\log \prod_{F\subseteq_0
D}\Big(p^{C_l}(i_F,i^*_{C_l\setminus F}\Big)^{(-1)^{|D\setminus
F|}},\;\;D\subseteq C_l,\;\;i_D\in {\cal I}^*_D\;.
\end{equation}
\end{definition}
\begin{definition}
\label{def:thetacond} For each residual $R_l$, $l=2,\ldots,k$,
 and
fixed $i_{S_l}\in {\cal I}_{S_l}$,
we define
\begin{equation}
\label{thetacond0}
 \theta^{R_l|i_{S_l}}(i_D)=\log
\prod_{F\subseteq_0 D}\Big(p^{R_l|i_{S_l}}(i_F, i^*_{R_l\setminus
F})\Big)^{(-1)^{|D\setminus F|}},\;\;D\subseteq R_l,\;\;i_D\in {\cal
I}^*_D\;.
\end{equation}
\end{definition}
Note that both $\theta^{C_l}(i_D)$ and  $\theta^{R_l|i_{S_l}}(i_D)$ are \lq \lq marginal\rq \rq{}
parameters, in the sense that they are
functions of probabilities  in the $C_l$-marginal table.

For any $A\subseteq V, B\subseteq V, B\cap A=\emptyset$ and any fixed $i_B\in {\cal I}_B$,
we also introduce the notation
\begin{eqnarray}
\theta({\cal D}^{C_1})&=&(\theta^{C_1}(i_D),\; D \subseteq C_1, \; i_D\in {\cal I}_D^*), \label{thetacond1} \\
n({\cal D}^{C_1})&=&(n(i_D),\; D \subseteq C_1, \; i_D\in {\cal I}_D^*),\nonumber
\end{eqnarray}
representing the log-linear parameters and cell-counts for the clique-$C_1$-table. Furthermore we will use
\begin{eqnarray}
\theta(i_B, {\cal D}^A)&=&\big(\theta^{A|i_{B}}(i_D),\;D\subseteq A, \;i_D\in {\cal I}_D^*\big)\label{thetacond+},\\
n(i_B, {\cal D}^A)&=&\big(n(i_B, i_D),\;D\subseteq A,\; i_D\in {\cal I}^*_D\big)\nonumber,
\end{eqnarray}
to represent the log-linear parameters  and the cell-counts respectively in the $i_B$-slice of the $A$-table.

We collect together the elements of (\ref{thetacond1}) and (\ref{thetacond+}) in a single
parameter that we call $\theta^{cond}$
\begin{equation}
\label{thetacond}
\theta^{cond}=\big(\theta({\cal D}^{C_1}),\;\;\;\theta(i_{S_l},{\cal D}^{R_l}),\;i_{s_l}\in {\cal
I}_{S_l},l=2,\ldots,k\big).
\end{equation}
Correspondingly we define the following canonical statistics
\begin{equation}
\label{statcond}
n^{cond}=\big(n({\cal D}^{C_1}),\;\;\;n(i_{S_l},{\cal D}^{R_l}),\;i_{s_l}\in {\cal
I}_{S_l},l=2,\ldots,k\big).
\end{equation}

Since $C_1$ and $R_l, l=2,\ldots,k$ are complete, we can apply
Proposition \ref{prop:basicmult} to each of the $C_1$-marginal
and $R_l$-conditional multinomials  in the $i_{S_l}$-slice of (\ref{cliquesepdecomp}).
We have the following lemma as an
immediate consequence of Proposition \ref{prop:basicmult}.
\begin{lemma}
\label{split}
The NEF-representation, in terms of the parametrization $\theta^{cond}$,
\begin{itemize}
\item
 of the marginal $C_1$-model
is given by
\begin{eqnarray}
\prod_{i_{C_1}\in {\cal I}_{C_1}}(p^{C_1}(i_{C_1}))^{n(i_{C_1})} &=&
\exp \{\<\theta({\cal D}^{C_1}),n({\cal D}^{C_1})\>-N \;k(\theta({\cal D}^{C_1}))\}.
\label{likthetac1}
\end{eqnarray}
\item
of the conditional $R_l$-model in the $i_{S_l}$-slice
%
%
is given by
\begin{eqnarray}
\hspace{-1mm}\prod_{i_{R_l}\in {\cal
I}_{R_l}}\Big(p^{R_l|i_{S_l}}(i_{R_l})\Big)^{n(i_{C_l})}&=&
 \exp \{\<\theta(i_{S_l},{\cal D}^{R_l}),n(i_{S_l},{\cal D}^{R_l})\>\nonumber\\
 &&\hspace{3cm}-n(i_{S_l})\;k(\theta(i_{S_l},{\cal D}^{R_l}))\}.
\label{liktheta}
\end{eqnarray}
\end{itemize}
\end{lemma}

%

Note that the number of parameters in
$\theta^{mod}$ and $\theta^{cond}$
is of course the same. Indeed each element of each one of the two
parametrizations is indexed by $i_D, D\in {\cal D}, i_D\in {\cal
I}_D^*$ either directly as for $\theta^{mod}$,  or through the
components $i_F,F\subseteq S_l$, $i_F\in {\cal I}^*_F$ and $i_D,
D\subseteq R_l, i_D\in {\cal I}^*_D$ as for
$\theta^{cond}$.
\vspace{2mm}

Since the clique marginal generalized log-odds ratios are also of
interest, we are now going to define a third parametrization of the
multinomial model in terms of the generalized log-odds ratios in (\ref{thetamarg}).
 Any marginal cell $i_{S_l}$ can be written as
$$i_{S_l}=(i_F,i^*_{S_l\setminus F})$$
where $F\in {\cal E}^*_{i_{S_l}}$. Accordingly, we define
\begin{eqnarray}
\theta({\cal D}_0^{S_l},{\cal D}^{R_l})&=&\big(\theta^{C_l}(i_F,i_D),\;D\subseteq
R_l,\;i_D\in {\cal I}^*_D,\;F\subseteq_0S_l,\;i_{F}\in {\cal I}^*_F\big)\nonumber\\
n({\cal D}_0^{S_l},{\cal D}^{R_l})&=&\big(n(i_F,i_D),\;D\subseteq
R_l,\;i_D\in {\cal I}^*_D,\;F\subseteq_0S_l,\;i_{F}\in {\cal I}^*_F\big)\nonumber
\end{eqnarray}
and
\begin{eqnarray}
\label{thetacliques}
\theta^{cliq}&=&\big(\theta({\cal D}^{C_1}),\;\;\;\theta({\cal D}_0^{S_l},{\cal D}^{R_l}),\;\;
\;l=2,\ldots,k\big)\;.\\
\label{ncliques}
n^{cliq}&=&\big(n({\cal D}^{C_1}),\;\;\;n({\cal D}_0^{S_l},{\cal D}^{R_l}),\;\;
\;l=2,\ldots,k\big)\;.
\end{eqnarray}
We note that for $F=\emptyset$,
$\theta^{C_l}(i_F,i_D)=\theta^{C_l}(i_D)$ and $n(i_F,i_D)=n(i_D)$. Clearly the number of
parameters in $\theta^{cliq}$ is the same as in $\theta^{cond}$.

The expression of the density in terms of this new parametrization
will be given in the next section, after we have derived the
relationship between the three parametrizations (\ref{thetamodel}),
(\ref{thetacond}) and (\ref{thetacliques}).

\subsection{Relationship between the various $\theta$ parametrizations}
The relationship between the three $\theta$ parametrizations is
given in the following proposition. To state the results succinctly,
let us also define, for any $F\subseteq V$ and $i_F\in {\cal I}^*_F$,
$$i_{\subseteq_0 F}=\{i_G, \;\;G\subseteq_0 F\}\;.$$
Then  for given $F\subseteq V$ and $i_F\in {\cal I}^*_F$,
 and $A\subseteq V$ such that $F\cap A=\emptyset$, we also define

\begin{equation}
\theta(i_{\subseteq_0F},{\cal D}^{A})=\big(\theta(i_G,
j_{L}),\;\;G\subseteq_0 F,\;L\subseteq A, \;j_L\in {\cal I}^*_L\big)
\end{equation}
and
\begin{equation}
\label{ked}
k(\theta(i_{\subseteq_0F}, {\cal D}^A))=\log \big(1+\sum_{L\subseteq A}
\exp \sum_{K\subseteq_0F\atop{H\subseteq L, \atop j_H\in {\cal I}^*_H}}\theta(i_K,j_H)\big).
\end{equation}
We note that for any $l=2,\ldots,k,$ and  $F\subseteq S_l$,
$$\theta(i_{\subseteq_0F},{\cal D}^{R_l})\subset\theta({\cal D}_0^{S_l},{\cal D}^{R_l}).
$$
\begin{prop}
Let $i_{D}\in {\cal I}^*_D$ and $D\subseteq C_l, D\cap
R_l\not =\emptyset$. Then

a) the relationship between $\theta^{cliq}$ and
$\theta^{cond}$ is
\begin{equation}
\label{relation} \theta^{C_l}(i_{D})=\sum_{F\subseteq_0 D\cap S_l}
(-1)^{|(D\cap S_l)\setminus
F|}\;\;\theta^{R_l|(i_F,i^*_{S_l\setminus F})}(i_{D\cap
R_l})
\end{equation}
which,  for $D\subseteq R_l$, is equivalent to
\begin{equation}
\label{relationinv} \theta^{R_l|(i_F,i^*_{S_l\setminus
F})}(i_{D})=\sum_{G\subseteq_0 F}\theta^{C_l}(i_G, i_{D})
\end{equation}

b) The relationship between $\theta^{cliq}$ and $\theta^{mod}$ is as follows.
Let $\{>l\}$ denote the set of $j\in \{l+1,\ldots,k  \}$ such that $C_l\cap C_j\not =\emptyset$.

(i) For $D\not \subseteq S_j$, for some $ j\in \{>l\}$,
\begin{equation}
\label{downwards0}
\theta(i_D)=\theta^{C_l}(i_D).
\end{equation}
(ii) For $D\subseteq S_m, m\in \{>l\}$
\begin{eqnarray}
\theta(i_D)&=&\theta^{C_l}(i_D)-\sum_{F\subseteq_0 D}(-1)^{|D\setminus F|}
k(\theta(i_{\subseteq_0F},{\cal D}^{C_{>l}}))\label{downwards}
\end{eqnarray}
where $C_{>l}=\cup_{m>l}(C_m\setminus C_l)$ and $k(\theta(i_{\subseteq_0F},{\cal D}^{C_{>l}}))$ is defined as in (\ref{ked}).

 Moreover, all $\theta(i_H,j_G)\in \theta(i_{\subseteq_0F},{\cal D}^{C_{>l}})$ are such
 that $H\cup G\subseteq C_m$ for some $m\in \{>l\}$ and is therefore either equal to
 $\theta^{C_m}(i_H,j_G)$ or can be expressed in terms of
 $\theta^{C_m}(i_E), m\in \{>l\}, E\subseteq C_m, i_E\in {\cal I}^*_E$.

\end{prop}
The proofs of (\ref{relation}) and (\ref{relationinv}) can easily be derived from Definitions \ref{def:thetamarg} and \ref{def:thetacond}. The proof of (\ref{downwards}), though, is not immediate and  is interesting. It is given in the appendix.


\vspace{0.3cm}
\noindent {\bf Remarks.}
\vspace{0.1cm}
1. Expression (\ref{relationinv}) is a generalization of the
relationship between conditional and marginal log-odds ratios for a
three way table given in Agresti (2002, p. 322).

2. According to (\ref{downwards}), $\theta(i_D)$ is a function of  $ \theta^{C_m}(j_H)$
such that $H\subseteq C_m$ for $m\geq l$ only.
This is going to be a crucial fact when we derive the reference prior of $\theta^{model}$
from the reference prior on $\theta^{cond}$ in the next section.

\vspace{2mm}

Relation (\ref{downwards}) is crucial for the derivation of the reference prior for $\theta^{cliq}$
in the next section and we therefore illustrate it here with an example.

\begin{example}
Consider a decomposable graphical model with the following perfect order of the cliques
$$C_1=\{a,b,c\},\;C_2=\{b,c,d\},\; C_3=\{c,d,e\}, C_4=\{e,f\},$$
having separators
$$S_2=\{b,c\},\;S_3=\{c,d\},\;S_4=\{e\}.$$
To simplify matters, let us assume the data are binary.
In this case we can simplify the notation since, because of the corner constraint conditions (see end of \S 2.2),
${\cal I}^*_D$ contains only one element for each $D$. Thus $\theta(i_D)$ can  more simply be  written $\theta(D)$.
Let us take $D=\{c,d\}$. We see that $D\subseteq C_2$ and $D\cap R_2=\{d\}\not =\emptyset$.
 Moreover  $C_{>2}=\{e,f\}$ and the set of $L\subseteq C_{>2}$ is equal to $\{e,f, ef\}$.
 Then according to (\ref{downwards}), it follows that
\begin{eqnarray}
\theta(cd)&=&\theta^{C_2}(cd)-\log (1+e^{\theta(e)+\theta(ec)+\theta(ed)+\theta(ecd)}+e^{\theta(f)}+e^{\theta(e)+\theta(ec)+\theta(ed)+\theta(ecd)+\theta(f)+\theta(ef)})\nonumber\\
&&+\log (1+e^{\theta(e)+\theta(ed)}+e^{\theta(f)}+e^{\theta(e)+\theta(ed)+\theta(f)+\theta(ef)})\nonumber\\
&&+\log (1+e^{\theta(e)+\theta(ec)}+e^{\theta(f)}+e^{\theta(e)+\theta(ec)+\theta(f)+\theta(ef)})\nonumber\\
&&-\log (1+e^{\theta(e)}+e^{\theta(f)}+e^{\theta(e)+\theta(f)+\theta(ef)})\nonumber
\end{eqnarray}
Since $$\theta(ec)=\theta^{C_3}(ec),\;\theta(ed)=\theta^{C_3}(ed),\;\theta(ecd)=\theta^{C_3}(ecd),\;\theta(ef)=\theta^{C_4}(ef),\;\theta(f)=\theta^{C_4}(f),$$
and according to (\ref{downwards}) again,
$$\theta(e)=\theta^{C_3}(e)+\log (1+e^{\theta(f)})-\log (1+e^{\theta(f)+\theta(ef)})=\theta^{C_3}(e)+\log (1+e^{\theta^{C_4}(f)})-\log (1+e^{\theta^{C_4}(f)+\theta^{C_4}(ef)}),$$
we see that $\theta(cd)$ can be expressed  in terms of $\theta^{C_m}(E), m\geq 2, E\subseteq C_m$.
\end{example}

 We will now give the
expression of the multinomial Markov model with respect to
$\theta^{cliq}$, using relation (\ref{relationinv}).
\begin{lemma}
Let $G$ be a
 decomposable
graph  with its cliques $C_1,\ldots,C_k$ arranged in a perfect
order. The NEF-representation of the multinomial Markov model in terms of the $\theta^{cliq}$
parametrization is given by
\begin{eqnarray}
\prod_{i\in{\cal I}}p(i)^{n(i)}&=&\exp \{\<\theta({\cal D}^{C_1}),n({\cal D}^{C_1})\>-N \;k(\theta({\cal D}^{C_1}))\nonumber\\
&\times& \prod_{l=2}^k\exp \{\<\theta({\cal D}_0^{S_l},{\cal D}^{R_l}),n({\cal D}_0^{S_l},{\cal D}^{R_l})\>\nonumber\\
&&\hspace{-2cm}-\sum_{F\subseteq_0S_l}\sum_{j_F\in {\cal I}^*_F}n(j_F)
\sum_{H\subseteq_0F}(-1)^{|F\setminus H|}k(\theta(j_{\subseteq_0H},{\cal D}^{R_l})) \}\label{multincliques}
\end{eqnarray}
\end{lemma}

From (\ref{multincliques}),  it appears that
under the multinomial Markov model,  the joint distribution of $n^{cliq}$  admits
 a conditional
reducibility structure, see Consonni and Veronese (2001); specifically,  it factorizes into the product
of $k$ conditional exponential families (save for the first term which is a marginal distribution),
in a recursive fashion according to the clique ordering.

\section{Reference priors}
\label{sec:reference priors}

In this section we shall derive reference priors for the various
parametrizations introduced in section 2. We shall only provide an
outline of the proofs of the derivation of our reference priors
since they follow the steps described in \S 2 and \S 4.2.1 of
Consonni et al. (2004).
 An important point to keep in mind is that a reference prior for a multidimensional parameter
depends on the grouping of its components, as well as the ordering of its groups:
specifically we order groups  according to inferential importance, while parameter-components that
 belong to the same
group are treated in a symmetric fashion.  For the   parametrizations considered in this paper, order
 will not matter, and thus the reference prior will only depend on the grouping-structure.

For a given graph $G$, let $C_1,\ldots,C_k$ represent a perfect ordering of the cliques.
We will first consider the reference prior for the collection of conditional
probabilities (including the marginal probabilities for clique $C_1$), $p^{cond}$ as in (\ref{cliquesepdecomp})
\begin{equation}
\label{pcond}
p^{cond}=(p^{C_1},p^{R_l|i_{S_l}}, \;i_{S_l}\in {\cal I}_{S_l},\;l=2,\ldots,k)\;,
\end{equation}
where
\begin{eqnarray}
p^{C_1}&=&\Big(p^{C_1}(i_{C_1}),\;\; i_{C_1}\in {\cal I}_{C_1}\Big) \label{pC1} \\
p^{R_l|i_{S_l}}&=&\Big(p^{R_l|i_{S_l}}(i_{R_l}),\;\;i_{R_l}\in {\cal I}_{R_l}\Big), \;
\label{pRlgiveniSl}\;
\end{eqnarray}
represent the collection of groups. Note that there are
$1+\sum_{l=2}^k  |{\cal I}_{S_l}|$
groups.
We remark that the nature of the parametrization $p^{cond}$
depends on the specific choice of the perfect numbering of the cliques
$C_1,\ldots,C_k$.
%

Next we will consider the reference priors for $\theta^{cond}, \theta^{cliq}, \theta^{mod}$
following  a parallel grouping-structure. We shall see that all these reference priors are strictly related,
so that a unified expression for all of them is possible.

\begin{prop}
\label{prop:refpriorp}
The reference prior for $p^{cond}$ relative to the grouping defined in (\ref{pcond}) is
\begin{equation}
\label{refp}
\pi^R_{p^{cond}}(p^{cond})\propto\Big(\prod_{i_{C_1}\in {\cal I}_{C_1}}p(i_{C_1})\Big)^{-\frac{1}{2}}\prod_{l=2}^k\prod_{i_{S_l}\in {\cal I}_{S_l}}\Big(\prod_{i_{R_l}\in {\cal I}_{R_l}}
p^{R_l|i_{S_l}}(i_{R_l})\Big)^{-\frac{1}{2}}.
\end{equation}
\end{prop}

  We note that the
reference prior is a product of Jeffreys' priors, one for each of the groups of $p^{cond}$.

\begin{pff}
In our setting, we simply need to derive the (Fisher) information matrix.
From (\ref{cliquesepdecomp}) it appears
that the likelihood function factorizes
into the product of terms, each involving exactly one group of $p^{cond}$;
furthermore each term is a saturated multinomial.
Accordingly
the information matrix is block-diagonal, and the determinant of each block,
using classic results,
%
is easily available. Specifically
the first one, corresponding to clique $C_1$, is given by
\begin{eqnarray}
N \Big(\prod_{i_{C_1}\in {\cal I}_{C_1}}p(i_{C_1})\Big)^{-1},
\end{eqnarray}
while for the remaining blocks the determinant is
\begin{eqnarray}
\label{eq:infoconditionalmultinomial}
E(n(i_{S_l})|p)\Big(\prod_{i_{R_l}\in {\cal I}_{R_l}}
p^{R_l|i_{S_l}}(i_{R_l})\Big)^{-1},\;\;\;i_{S_l}\in {\cal I}_{S_l}, l=2,\ldots,k\;.
\end{eqnarray}
 Because of the perfect ordering the cliques,
 $S_l\subseteq C_j$ for some $j<l$,
so that the expected value $E(n(i_{S_l})|p)$ is a function of  parameters
only belonging to groups preceding the $l$-th one.
%
%

Following the theory summarized in Consonni et al. (2004, sect. 2), the
reference prior is given by the
square root of the
product of the block-determinants,
excluding the terms $E(n(i_{S_l})|p)$,  and
the  result is established.
\end{pff}

We now emphasize three properties
 of the reference prior for $p^{cond}$.
First of all, since the information matrix is block-diagonal,  the reference prior is order-invariant, i.e.
it does not depend
on the order of   the groups.
Secondly, we  remark that there exists also some degree of invariance with respect to grouping. Specifically, if
if we lumped together
%
%
in one single block
 all the $i_{S_l}$ terms  $p^{R_l|i_{S_l}}
 ,\;\;i_{S_l}\in {\cal I}_{S_l}$, the reference prior would not change.
 This feature will turn out to be useful later on when deriving reference priors for alternative parametrizations.

Third, we  remark that the distribution $\pi^R_{p^{cond}}$
belongs to a family conjugate to the likelihood for $p^{cond}$, see   (\ref{cliquesepdecomp}).
Accordingly its hyper-parameters can be interpreted in terms of \lq \lq prior counts\rq \rq{}; the latter
however cannot be recovered as the margins of an fictitious overall table. Indeed, each cell
in the $C_1$-table, as well as in the $i_{S_l}$ slice of the $R_l$-table, has a prior count equal to $1/2$, irrespective of the dimension of the subtables and of the overall table.
Finally, the prior is proper, since it is a product of Dirichlet priors, one for each block, each Dirichlet
being indexed by a vector of hyper-parameters with entries all equal to $1/2$.

\vspace{.5mm}

We now turn to the derivation of the reference priors for the three $\theta$ parametrizations described in
\S \ref{sec:parametrizations}.
Central to our arguments below is the following basic fact about reference priors that we
shall use repeatedly. Let $\lambda$ be a parameter grouped into components
$\lambda=(\lambda_1,\ldots,\lambda_k)$, where
$\lambda_i$ is typically  a vector. We assume that the above groups
are arranged in increasing order of inferential importance.   Let $\phi=(\phi_1,\ldots,\phi_k)$ be a reparametrization,
i.e. a one-to-one function of $\lambda$ with
$\phi_i$ having the same dimension as $\lambda_i$.
Suppose that, for each $l=1,\ldots,k$, $\phi_l=h_l(\lambda_1,\ldots,\lambda_l)$,
for some function $h_l$: we say in this case that the map $\lambda \mapsto \phi$
is block-lower triangular. Then the reference prior for $\phi$ can be obtained  from that
of  $\lambda$ simply by a change of variable. For details and references see again
Consonni et al. (2004, \S 4.2.1)
%
%
%
An important special case occurs when $\phi_l=h_l(\lambda_l)$: in this case
we say that the map is block-wise one-to-one.

We start by expressing  $p^{cond}$ in terms of the $\theta^{cond}$. In order to achieve this goal,
%
it is  convenient to define the parameters
\begin{eqnarray}
\xi^{C_1}(i_D)&=&\sum_{F\subseteq D}\theta^{C_1}(i_F),\;\;i_D\in {\cal I}^*_D \label{xiC1} \\
\xi^{R_l|i_F}(i_F,i_D)&=&\sum_{L\subseteq D}\theta^{R_l|(i_F,i^*_{S_l\setminus F})}(i_L)\label{xitothetacond}\\
&=&\sum_{L\subseteq D}\sum_{H\subseteq_0 F}\theta^{C_l}(i_H,i_L),\;\;i_F\in {\cal I}^*_F, i_L\in {\cal I}^*_L,\;
D\subseteq R_l\label{xitothetacliques}
\end{eqnarray}
We let
\begin{equation}
\xi=\big(\xi^{C_1},\;\;\xi^{R_l|(i_F,i^*_{S_l\setminus F})},\;F\subseteq S_l, i_F\in {\cal I}^*_F,\;\;l=2,\ldots,k\big)
\end{equation}
where
\begin{eqnarray}
\xi^{C_1}&=&(\xi^{C_1}(i_D),\;D\subseteq C_1)\nonumber\\
\xi^{R_l|(i_F,i^*_{S_l\setminus F})}&=&(\xi^{R_l|(i_F,i^*_{S_l\setminus F})}(i_F,i_D),\;D\subseteq R_l, i_D\in {\cal I}^*_D)\;,\nonumber
\end{eqnarray}
The mapping between $p^{cond}$  and $\xi$ is block-wise one-to-one.
As a consequence the reference prior on $\xi$ can be deduced from that of $p^{cond}$ as

\begin{eqnarray}
\pi^R_{\xi}(\xi)=\pi^R_{p^{cond}}(p^{cond}(\xi))|J_{p^{cond}}(\xi)|
\end{eqnarray}
where $J_{p^{cond}}(\xi)$ is the Jacobian of the transformation $p^{cond} \mapsto \xi$. It can be verified that
\begin{equation}
\det \left(\frac{d p^{cond}}{d \xi}\right)=\Big(\prod_{i_{C_1}\in {\cal I}_{C_1}}p(i_{C_1})(\xi^{C_1})\Big)\prod_{l=2}^k\prod_{i_{S_l}\in {\cal I}_{S_l}}\Big(\prod_{i_{R_l}\in {\cal I}_{R_l}}
p^{R_l|i_{S_l}}(i_{R_l})(\xi^{R_l|i_{S_l}})\Big)
\end{equation}
so that the induced reference prior for $\xi$ is
\begin{equation}
\pi^R_{\xi}(\xi) \propto \Big(\prod_{i_{C_1}\in {\cal I}_{C_1}}p(i_{C_1})(\xi^{C_1})\Big)^{-\frac{1}{2}+1}\prod_{l=2}^k\prod_{i_{S_l}\in {\cal I}_{S_l}}\Big(\prod_{i_{R_l}\in {\cal I}_{R_l}}
p^{R_l|i_{S_l}}(i_{R_l})(\xi^{R_l|i_{S_l}})\Big)^{-\frac{1}{2}+1}
\end{equation}
We shall also need the following result
%
%
%
which can be easily derived from Definitions (\ref{def:thetamarg}) and (\ref{def:thetacond}) and Moebius inversion formula.
\begin{lemma}
For $i_{C_1}=(i_F,i^*_{C_1\setminus F})$,
\begin{equation}
\label{pigeneralgcl}
p^{C_l}(i_{C_l})=\frac{\exp \xi^{C_l}(i_F)}{1+\sum_{H\subseteq C_1}\sum_{j_H\in {\cal I}^*_H}\exp\xi^{C_l}_D(j_H)}
\;,
\end{equation}
For $i_{S_l}$ and $i_{R_l}=(i_G, i^*_{R_l\setminus G})$ given,
\begin{equation}
\label{pigeneralgrl}
p^{R_l|i_{S_l}}(i_{R_l})=\frac{\exp \xi^{R_l|i_{S_l}}(i_G)}{1+
\sum_{H\subseteq R_l}\sum_{j_H\in {\cal I}^*_H}\exp \xi^{R_l|i_{S_l}}_D(j_H)}
\;,
\end{equation}
As particular cases, we have
\begin{equation}
\label{piemptygcl}
p^{C_l}(i^*_{C_l})=\frac{1}{1+\sum_{H\subseteq C_1}\sum_{j_H\in {\cal I}^*_H}\exp\xi^{C_l}_D(j_H)}
\;,\;\;
\end{equation}
and
\begin{equation}
\label{piemptygrl}
p^{R_l|i_{S_l}}(i^*_{R_l})=\frac{1}{1+
\sum_{H\subseteq R_l}\sum_{j_H\in {\cal I}^*_H}\exp \xi^{R_l|i_{S_l}}_D(j_H)}\;\;
\end{equation}
\end{lemma}

%
Since the reference
priors of the three $\theta$-parametrizations   are
structurally equivalent we shall provide the result in a unified statement.
\begin{theorem}
\label{theo:refvarioustheta}
The reference prior for
\begin{itemize}
\item [a)]
$\theta^{cond}$,
relative to the grouping defined in (\ref{thetacond})
\item [b)]
$\theta^{cliq}$, relative
to the grouping defined in (\ref{thetacliques})
\item [c)]
$\theta^{mod}$, relative to the following grouping
\begin{eqnarray}
\label{order}
&&\hspace{-3cm}\tilde{\theta}^{C_1}=(\theta(i_D), D\subseteq C_1, i_D\in {\cal I}^*_D),
 \tilde{\theta}^{C_l}=(\theta(i_D), D\subseteq C_l, D\cap R_l\not =\emptyset), l=2,\ldots,k\;.
\end{eqnarray}
\end{itemize}
is proportional to
\begin{equation}
\label{refvarioustheta}
\Big(\prod_{i_{C_1}\in {\cal I}_{C_1}}p(i_{C_1})(\cdot)\Big)^{\frac{1}{2}}
\prod_{l=2}^k\prod_{i_{S_l}\in {\cal I}_{S_l}}\Big(\prod_{i_{R_l}\in {\cal I}_{R_l}}
p^{R_l|i_{S_l}}(i_{R_l})(\cdot)\Big)^{\frac{1}{2}},
\end{equation}
where the probabilities $p(i_{C_1})(\cdot)$ and $p^{R_l|i_{S_l}}(i_{R_l})(\cdot)$ are
understood to be expressed  in terms of the relevant $\theta$-parametrization, using
(\ref{pigeneralgcl})-(\ref{piemptygrl}) together with
i)(\ref{xiC1})-(\ref{xitothetacond})
for $\theta^{cond}$;
ii) (\ref{xitothetacliques}) for $\theta^{cliq}$.
iii)  (\ref{relationinv}), (\ref{downwards0})  and (\ref{downwards}) for
$\theta^{mod}$.

More explicitly, the reference prior
\begin{itemize}
\item
for $\theta^{cond}$ is given by the product of (\ref{likthetac1})
and (\ref{liktheta}),
with the understanding that
the  counts in these formulas are replaced by fictitious prior counts which we write as
$\tilde{n}(i_D), \tilde{N}$ and so on. More precisely, we have
\begin{eqnarray*}
 \tilde{n}(i_D)=\frac{|{\cal I}_{C_1\setminus D}|}{2},\;\;\tilde{N}=\frac{|{\cal I}_{C_1}|}{2},\;\;
\end{eqnarray*}
and
\begin{eqnarray*}
\tilde{n}(i_{S_l},i_D)=\frac{|{\cal I}_{R_l\setminus D}|}{2},\;\;\;
\tilde{n}(i_{S_l})=\frac{|{\cal I}_{R_l}|}{2}.
\end{eqnarray*}
\item
for $\theta^{cliq}$ is given by
(\ref{multincliques})
where for $l=1$
\begin{eqnarray*}
 \tilde{n}(i_D)=\frac{|{\cal I}_{C_1\setminus D}|}{2}\;\;\;\mbox{and}\;\;\;\tilde{N}=\frac{|{\cal I}_{C_1}|}{2},\;\;
\end{eqnarray*}
and for $l=2,\ldots,k$
\begin{eqnarray*}
\tilde{n}(j_F,i_D)=\frac{|{\cal I}_{S_l\setminus F}||{\cal I}_{R_l\setminus D}|}{2}\;\;\;\mbox{and}\;\;\;
\tilde{n}(j_F)=\frac{|{\cal I}_{S_l\setminus F}||{\cal I}_{R_l}|}{2}
\end{eqnarray*}
\item
for $\theta^{mod}$ can be obtained from  that of $\theta^{cliq}$ above
by expressing it in terms of $\theta({\cal D})$
using (\ref{downwards0}) and (\ref{downwards}).
\end{itemize}
\end{theorem}

\begin{pff}
a)
Because of (\ref{xiC1}) and (\ref{xitothetacond}) it is immediate to verify that the map $\xi \mapsto \theta^{cond}$
is block-wise one-to-one; moreover the Jacobian is equal to one. Accordingly the reference prior
for $\theta^{cond}$ will be exactly as that for $\xi$, with the only difference that the probabilities involved
will be expressed as functions of $\theta^{cond}$.

b) Similarly to what happened for the reference prior for $p^{cond}$,
the reference prior
for $\theta^{cond}$ is unchanged if, for each $l=2,\ldots,k$, we lump together the groups
labeled by $i_{S_l}\in {\cal I}_{S_l}$, and thus only regard $\theta^{cond}$ as made up of $k$ groups.
In this way
the transformation from $\theta^{cond}$ to $\theta^{cliq}$
is block-wise one-to-one, and thus the reference prior for $\theta^{cliq}$
is equal to that induced from the reference prior $\theta^{cond}$. Moreover, the transformation is
 linear so that the Jacobian is constant, and thus the result follows.

c) We see that  the groupings in (\ref{order}) are exactly parallel to those in $\theta^{cliq}$.
From (\ref{downwards}) we also see that the $l$-th  group in $\theta^{mod}$ is a function
of the subsequent $l,l+1,\ldots,k$ groups in $\theta^{cliq}$.
This
defines a block-upper triangular transformation, which can be turned into a block-lower triangular one
by reversing the order of the groups in $\theta^{cliq}$. Since  the
reference prior on $\theta^{cliq}$ is invariant to group-ordering,
we conclude that the reference prior
on $\theta^{mod}$ can be obtained from that of $\theta^{cliq}$ by a change-of-variable. Now notice
that the Jacobian is 1, again using (\ref{downwards}), so that the result is proved.
Finally, the expressions of the fictitious counts are derived by inspection.

\end{pff}

We remark that, similarly to what happened for $p^{cond}$, the reference prior for each of
the three $\theta$-parametrizations is also a conjugate prior
  and is proper, being the transformation of a proper prior on $p^{cond}$.

\section{Parametrizations and reference priors associated to  a cut}

The reference priors obtained in the previous section were based on a grouping of the parameters
defined by the
structure of the graph,  essentially through a perfect order  of the cliques (and consequently
of residuals and separators).

Now suppose we are interested in a particular subset $A \subseteq V$ of the variables, and that we would like
to consider a reference prior which groups together precisely the  parameters of the marginal
distribution referring to $A$.
%
%
We show in this section how this can be done if the Markov model  ${\cal M}_G$ is \emph{collapsible} onto $A$,
equivalently if $A$  represents  a \emph{cut} for the joint distribution.

Asmussen and Edwards (1983) consider the concept of collapsibility
for contingency tables. If the set of factors for the table are
indexed by $\gamma\in V$ and if $A\subseteq V$, we say that $G$ is
collapsible onto $A$ if the multinomial model ${\cal M}_{G_A}$,
 Markov with respect to the induced subgraph $G_A$,
  is the same as the model obtained by marginalising the given model ${\cal M}_G$, Markov with respect to $G$,
  over the $A$-table.

Frydenberg (1990, Theorem 5.4) has shown that the model for the
random vector $Y$, Markov with respect to $G$, is collapsible
onto $A$ if and only if  the sub-vector $Y_A$ is a cut (for simplicity we shall also say that $A$ is a
cut). Cuts in exponential families have been introduced in
Barndorff-Nielsen (1978) and studied in several further articles such
as Barndorff-Nielsen and Koudou (1995).

%
%
%

A very useful result, due to
Asmussen and Edwards (1983), is that $A$ will induce a cut
if and only if the boundary of every connected
component of $V\setminus A$ has a complete boundary in $G$.
%

The following lemma gives us the factorization of ${\cal M}_G$ with respect to the cut $A$ and the
connected components of $G_{V\setminus A}$.
\begin{lemma}
\label{lem:factorization_cut}
Let $A$ be a cut. Let $C_1',\ldots, C_q'$ be a perfect ordering of the cliques of $G_A$,
the graph induced by $A$.
Let $B_l,l=1,\ldots,p$ be the connected components of $G_{V\setminus A}$.
Let $C^{l}_j, j=1,\ldots,m_l$ be the cliques of the induced graph $G_{B_l\cup \partial B_l}, l=1,\ldots,p$.
The  multinomial model ${\cal M}_G$, Markov with respect to $G$, can be factorized as follows
\begin{eqnarray}
\hspace{-1cm}\prod_{i\in {\cal I}}p(i)^{n(i)}&=&\prod_{i_{C_1'}\in {\cal I}_{C_1'}}(p^{C_1'}(i_{C_1'}))^{n(i_{C_1'})}\prod_{l=1}^q
\prod_{i_{S_l'}\in {\cal I}_{S_l'}}\prod_{i_{R_l'}\in {\cal I}_{R_l'}}(p^{R_l'|i_{S_l'}}(i_{R_l'}))^{n(i_{C_l'})}\nonumber\\
\hspace{-8mm}\times  \prod_{l=1}^p && \hspace{-1cm} \left(\prod_{i_{\partial B_l}}\prod_{i_{C_1^{(l)}\setminus \partial B_l}}
(p^{C_1^{(l)}\setminus \partial B_l|i_{\partial B_l
}}(i_{C_1^{(l)}\setminus \partial B_l}))^{n(i_{C_1^{(l)}})}
\prod_{j=2}^{m_l}\prod_{i_{S_j^{(l)}}}\prod_{i_{R_j^{(l)}}}(p^{R_j^{(l)}|i_{S_j^{(l)}}}(i_{R_j^{(l)}}))^{n(i_{C_j^{(l)}})}
\right)\hspace{5mm}\label{factcuta}
\end{eqnarray}
\end{lemma}

\begin{pff}
For simplicity of exposition, some statements concerning the random variables
associated to a set, will be simply stated in terms of the set itself.
If $A$ is a cut,
A separates the connected components of $V\setminus A$;
by Theorem 2.8 of Dawid and Lauritzen (1993), this implies that the $B_l$'s are
mutually conditionally independent  given $A$.
Moreover since $A$ is a cut, the boundary of $B_l$ is a complete subset of $A$ and, of course,
it separates $B_l$ from $V\setminus (B_l\cup \partial B_l)$.
Therefore the overall multinomial Markov model
factorizes as the product of the $A$-marginal multinomial model, Markov
with respect to  ${\cal M}_{G_A}$,
and the product of the conditional multinomial distributions
of the $B_l$s  given $i_{\partial B_l}$,   $l=1,\ldots,p$,.


 Since the the marginal model for $A$
 is
 Markov with respect to the graph $G_A$,
 it factorizes according to a perfect order of the cliques of $G_A$, in
 parallel to what was done in \S \ref{sec:reference priors}: this proves the first line of
 (\ref{factcuta}).

Let us now consider the expression for the second line of  (\ref{factcuta}).
As recalled above, this is given by the product of the conditional multinomial models
for $B_l, l=1,\ldots,p$ given $i_{\partial B_l}$.
For any $l\in \{1,\ldots,p\}$, as a subgraph of $G$, the induced graph $G_{B_l\cup \partial B_l}$ is decomposable.
%
Moreover the marginal model for $B_l\cup \partial B_l$ is  Markov w.r.t. $G_{B_l\cup \partial B_l}$.
This happens because $B_l\cup \partial B_l$ is itself a cut, since the boundary of
each connected component of $G_{V \setminus (B_l\cup \partial B_l)}$ clearly belongs to $\partial B_l$
which is complete.
Therefore the marginal distribution ${\cal M}_{G_{B_l\cup
\partial B_l}}$ factorizes according to a perfect order of the
cliques of $G_{B_l\cup \partial B_l}$. Since $\partial B_l$ is
complete, it must belong to a clique $C_1^{(l)}$ of $G_{B_l\cup
\partial B_l}$ and by Proposition 2.29 of Lauritzen (1996), we know
that we can take this clique as the first in a perfect order
$C_i^{(l)}, i=1,\ldots, m_l$ of the cliques of $G_{B_l\cup \partial
B_l}$.

%
The marginal multinomial distribution ${\cal M}_{G_{B_l\cup \partial B_l}}$ can therefore be written as
\begin{eqnarray}
\prod_{j=1}^{m_l}\prod_{i_{C^{(l)}_j}}(p^{C_j^{(l)}}(i_{C^{(l)}_j}))^{n(i_{C^{(l)}_j})}&=&
\prod_{i_{C_1^{(l)}}}
p^{C_1^{(l)}}(i_{C_1^{(l)}})^{n(i_{C_1^{(l)}})}
\prod_{j=2}^{m_l}\prod_{i_{S_j^{(l)}}}\prod_{i_{R_j^{(l)}}}(p^{R_j^{(l)}|i_{S_j^{(l)}}}(i_{R_j^{(l)}}))^{n(i_{C_j^{(l)}})}\nonumber\\
&=&\prod_{i_{\partial B_l}}\left((p^{\partial B_l}(i_{\partial B_l}))^{n(i_{\partial B_l})}\prod_{i_{C^{(l)}_1\setminus \partial B_l}}(p^{C^{(l)}_1\setminus \partial B_l|{i_{\partial B_l}}}(i_{C^{(l)}_1\setminus \partial B_l}))^{n(i_{C^{(l)}_j})}\right)\nonumber\\
&\times
&\prod_{j=2}^{m_l}\prod_{i_{R_j^{(l)}}}\prod_{i_{S_j^{(l)}}}p^{R_j^{(l)}|i_{S_j^{(l)}}}(i_{R_j^{(l)}})^{n(i_{C_j^{(l)}})}\nonumber
\end{eqnarray}
and therefore
the model for $B_l$ conditional on  $i_{\partial B_l}$ is equal to
\begin{eqnarray}
\label{eq:modelBlgivenpartialBl}
&&\prod_{i_{C^{(l)}_1\setminus \partial B_l}}(p^{C^{(l)}_1\setminus \partial B_l|{i_{\partial B_l}}}(i_{C^{(l)}_1\setminus \partial B_l}))^{n(i_{C^{(l)}_j})}
\prod_{j=2}^{m_l}\prod_{i_{R_j^{(l)}}}
\prod_{i_{S_j^{(l)}}}p^{R_j^{(l)}|i_{S_j^{(l)}}}(i_{R_j^{(l)}})^{n(i_{C_j^{(l)}})}
.
\end{eqnarray}
Since this is true for all $B_l$, the result is established.
\end{pff}
\begin{example}
Suppose that the  joint distribution of the 11 variables numbered
consecutively from 1 to 11 is Markov with respect to the
decomposable graph $G$ as given in Fig. \ref{Fig1}.

\begin{figure}
\begin{center}
\includegraphics[width=12cm]{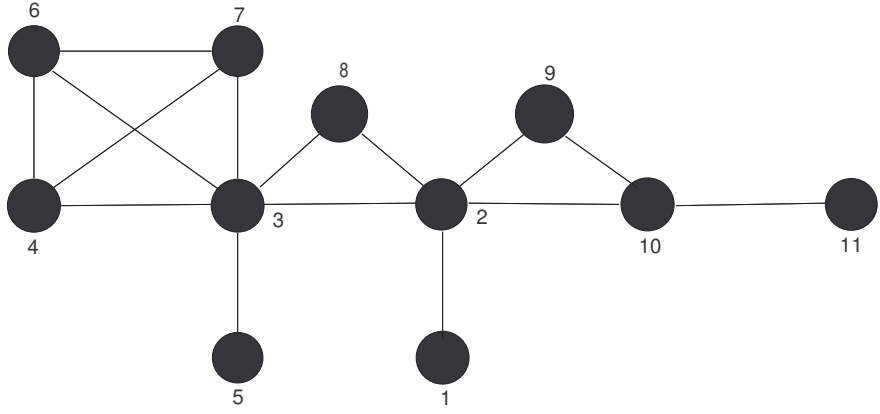}\\
\medskip
\parbox{15cm}
{
\caption[]{The decomposable graph for Example 4.1}.\\
\label{Fig1} }
\end{center}
\end{figure}

Consider the subset of
variables given by $A=\{1,2,3,4 \}$.
A perfect ordering of the cliques of the induced sub-graph $G_A$ is
\begin{eqnarray}
C^{\prime}_1=\{ 1,2 \}, \,C^{\prime}_2=\{ 2,3 \}, \,C^{\prime}_3=\{ 3,4 \},
\end{eqnarray}
so that
$S^{\prime}_2=\{ 2 \},\, \,S^{\prime}_3=\{ 3 \},
R^{\prime}_2=\{ 3 \}, \, \,R^{\prime}_3=\{ 4 \}$.
The connected components $B_l$ of $G_{V\setminus A}$, their boundary $\partial B_l$ together with
 the cliques $C_j^{l}$  of $G_{B_l \cup \partial B_l}$
are
\begin{center}
\begin{tabular}{c|c|c|c|c}
$l$ & $B_l$& $\partial B_l$ & ${B_l \cup \partial B_l}$ & $C_j^{(l)}$
\\ \hline
1& \{9,10,11\} & \{2\}& \{2,9,10,11\} & $C^{(1)}_1=\{2,9,10\},\,C^{(1)}_2=\{10,11\}$
\\
 2& \{8\} & \{2,3\}& \{2,3,8\} & $C^{(2)}_1=\{2,3,8\}$\\
3& \{5\} & \{3\}& \{3,5\} & $C^{(3)}_1=\{3,5\}$\\
4& \{6,7\} & \{3,4\}& \{3,4,6,7\} & $C^{(4)}_1=\{3,4,6,7\}$ \\
 \hline
\end{tabular}
\end{center}
\vspace{.2cm}
A graphical display of $G_A$ and its connected components is given in Fig. \ref{Fig2}.

\begin{figure}
\begin{center}
\includegraphics[width=12cm]{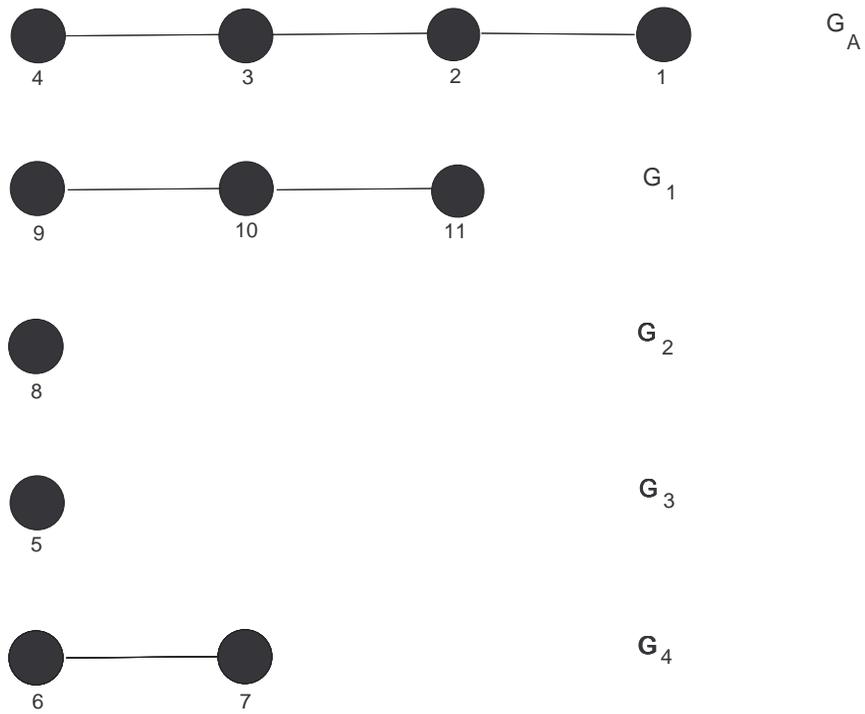}\\
\medskip
\parbox{15cm}
{
\caption[]{The decomposable graph $G_A$ associated to a cut $A$ and the connected
components of $G_{V \setminus A}$  for Example 4.1}.\\
\label{Fig2} }
\end{center}
\end{figure}

Accordingly,  the multinomial model, Markov with respect to $G$, can be
factorized using Lemma \ref{lem:factorization_cut} as
%

\begin{eqnarray*}
%
\prod_{i\in {\cal I}}p(i)^{n(i)}
&=& \\
  &&
\hspace{-2cm}
  \prod_{i_{C_1'}\in {\cal I}_{C_1'}}(p^{C_1'}(i_{C_1'}))^{n(i_{C_1'})}
\prod_{i_{S_2'}\in {\cal I}_{S_2'}}\prod_{i_{R_2'}\in {\cal I}_{R_2'}}
(p^{R_2'|i_{S_2'}}(i_{R_2'}))^{n(i_{C_2'})}
\prod_{i_{S_3'}\in {\cal I}_{S_3'}}\prod_{i_{R_3'}\in {\cal I}_{R_3'}}
(p^{R_3'|i_{S_3'}}(i_{R_3'}))^{n(i_{C_3'})}\\
&\times&
 \prod_{i_{\partial B_1}}\prod_{i_{C_1^{(1)}\setminus \partial B_1}}
(p^{C_1^{(1)}\setminus \partial B_1|i_{\partial B_1
}}(i_{C_1^{(1)}\setminus \partial B_1}))^{n(i_{C_1^{(1)}})}
\prod_{i_{S_2^{(1)}}}\prod_{i_{R_2^{(1)}}}(p^{R_2^{(1)}|i_{S_2^{(1)}}}(i_{R_2^{(1)}}))^{n(i_{C_2^{(1)}})}
\\
&\times& \prod_{i_{\partial B_2}}\prod_{i_{C_1^{(2)}\setminus \partial B_2}}
(p^{C_1^{(2)}\setminus \partial B_2|i_{\partial B_2
}}(i_{C_1^{(2)}\setminus \partial B_2}))^{n(i_{C_1^{(2)}})}
\\
&\times& \prod_{i_{\partial B_3}}\prod_{i_{C_1^{(3)}\setminus \partial B_3}}
(p^{C_1^{(3)}\setminus \partial B_3|i_{\partial B_3
}}(i_{C_1^{(3)}\setminus \partial B_3}))^{n(i_{C_1^{(3)}})}
\\
&\times& \prod_{i_{\partial B_4}}\prod_{i_{C_1^{(4)}\setminus \partial B_4}}
(p^{C_1^{(4)}\setminus \partial B_4|i_{\partial B_4
}}(i_{C_1^{(4)}\setminus \partial B_4}))^{n(i_{C_1^{(4)}})}
%
\end{eqnarray*}

\end{example}

We now provide the expression for the reference prior associated to a cut.

\begin{theorem}
Let $A$ be a cut
and consider
 the parametrization associated to  $A$
$$
p^{cut}_A= (p^{\prime,\,cond}_A,p^{\prime,\,cond}_{V \setminus A |A}),
$$
where
\begin{eqnarray}
p^{\prime,\,cond}_A & = &(p^{C_1'},p^{R_l'|i_{S_l'}}, l=1,\ldots, q, \,{i_{S_l'}\in {\cal I}_{S_l'}})
\label{eq:cut-parameter-1}\\
p^{\prime,\,cond}_{V \setminus A |A} & = &
(
p^{C_1^{(l)}\setminus \partial B_l|i_{\partial B_l}}, \,
 i_{\partial B_l} \in {\cal I}_{\partial B_l}; \,
p^{R_j^{(l)}|i_{S_j^{(l)}}}
l=1,\ldots,p,\nonumber\\
&&\hspace{6cm}
\,j=2, \ldots, m_l, \, i_{S_j^{(l)}} \in {\cal I}_{S_j^{(l)}}
),
\label{eq:cut-parameter-2}
\end{eqnarray}
using the notation presented in Lemma \ref{lem:factorization_cut}.
The reference prior for  $p^{cut}_A$, relative to the grouping (\ref{eq:cut-parameter-1}) and
(\ref{eq:cut-parameter-2}), is
\begin{eqnarray*}
\pi^R_{p^{cut}_A}(p^{cut}_A) &\propto&
\prod_{i_{C_1'}\in {\cal I}_{C_1'}}p^{C_1'}(i_{C_1'})^{-\frac{1}{2}}\prod_{l=1}^q
\prod_{i_{S_l'}\in {\cal I}_{S_l'}}\prod_{i_{R_l'}\in {\cal I}_{R_l'}}(p^{R_l'|i_{S_l'}}(i_{R_l'}))^{-\frac{1}{2}}
\\
 &\times& \prod_{l=1}^p\left(\prod_{i_{\partial B_l}}\prod_{i_{C_1^{(l)}\setminus \partial B_l}}
(p^{C_1^{(l)}\setminus \partial B_l|i_{\partial B_l
}}(i_{C_1^{(l)}\setminus \partial B_l}))^{^{-\frac{1}{2}}}
\prod_{j=2}^{m_l}\prod_{i_{S_j^{(l)}}}\prod_{i_{R_j^{(l)}}}(p^{R_j^{(l)}|i_{S_j^{(l)}}}(i_{R_j^{(l)}}))^{^{-\frac{1}{2}}}
\right)
\end{eqnarray*}

\end{theorem}
We emphasize that, also for this case, the prior admits a conjugate structure and is proper, being a product
of Jeffreys' priors.

\begin{pff}
Using Lemma \ref{lem:factorization_cut} the likelihood factorizes into a product of two general terms,
one related to the marginal distribution of $A$ indexed by $p^{\prime,\,cond}_A$,
the other related to the conditional distribution
of $V \setminus A$ given $A$ indexed by $p^{\prime,\,cond}_{V \setminus A |A}$. The two
groups of parameters are  variation and likelihood independent, so that the   information matrix is two-block-diagonal.
The marginal distribution related to $A$ is a $G_A$-Markov model, with $G_A$ decomposable, and therefore the
corresponding reference prior is exactly as in the general decomposable case of Proposition  \ref{prop:refpriorp}.
This yields the first line of the kernel of the reference prior
%

To prove the second line,
we have  to consider the second block of the information matrix.
This actually further decomposes into $p$ diagonal blocks, one for each
connected component $B_l$
. Consider  the block corresponding to
the model for  $B_l$ conditional on $\partial B_l, l=1,\ldots,p$ (see (\ref{eq:modelBlgivenpartialBl})).
Each block represents the information of a saturated conditional multinomial.
In particular the first term has cell-probabilities  $p^{C_1^{(l)}\setminus \partial B_l|i_{\partial B_l}}$ and
 $n(i_{\partial B_l})$ trials, while the remaining terms have cell-probabilities
 $p^{R_j^{(l)}|i_{S_j^{(l)}}}$ and  $n(i_{S_j^{(l)}})$ trials.
 The expression of the corresponding term in the information matrix will therefore be as
in the general  conditional saturated multinomial, see (\ref{eq:infoconditionalmultinomial}).
Finally, the expectation of  $n(i_{\partial B_l})$  depends only on the parameter
$p^{\prime,\,cond}_A$ since $\partial B_l \in A$, and
  similarly the expectation of $n(i_{S_j^{(l)}})$ does not depend on the parameter
   $p^{R_j^{(l)}|i_{S_j^{(l)}}}$
  specific to the component because of the perfect ordering of the cliques. Therefore, in both
  cases   the term corresponding to  the expectation factors out of the determinant and the
  proof is complete.
\end{pff}

\section{Discussion}
In this paper we have considered several alternative
parametrizations for discrete decomposable graphical models. First
of all we have described a parametrization   in terms of
 conditional cell-probabilities. Next we have derived three alternative representations in terms
 of natural exponential families,
whose canonical parameters
represent generalized log-odds ratios relative to suitable
cell-probabilities. Specifically, $\theta^{mod}$ refers to the
joint probabilities of the full table and has been previously
used, see e.g. Dellaportas and Forster (1999), Dellaportas and
Tarantola (2005) and Liu and Massam (2007) but we think that our derivation and
interpretation makes its interpretation clearer. The parametrizations
$\theta^{cond}$ and $\theta^{clique}$, on the other hand, refer to
marginal  sub-tables and are quite distinct from those
traditionally employed in graphical log-linear modelling. Indeed
they are rather related to the concept of marginal models,
 see e.g. Bergsma and Rudas (2002), Lang and Agresti (1994) and
Glonek and McCullagh (1995).

A reference prior for each of the above parametrizations was
constructed. In particular the prior for the conditional
cell-probabilities of the residuals given the separators is a product,
of Jeffreys' priors.
 We showed that all reference priors are
coherent, i.e. each is equivalent to any other one. This happens because the
grouping structure is such that the transformation between any two
parametrizations is either block-diagonal or block-lower
triangular.
A notable feature is that all  reference priors are proper.
Another property is that they belong to a conjugate
family, which facilitates prior-to-posterior updating.
  The conjugacy feature  is consistent with
previous  results, see Consonni et al. (2004), wherein reference priors for suitable parametrization of
NEFs having a simple quadratic variance (such as the saturated
multinomial) were derived and shown to belong to (enriched)
conjugate families. Our paper shows that this result continues to
hold  also for multinomial decomposable models, whose variance
function is not quadratic. With hindsight, this is not surprising,
because of the recursive factorization into products of
conditional saturated multinomial models that holds when $G$ is decomposable.
We have also considered a parametrization, and the corresponding reference prior, associated
to a cut. This can be especially useful whenever interest focuses on the parameters of a marginal
table, e.g. because of their inferential interest.

Throughout the paper we assumed
 a given graphical model, and  constructed reference priors essentially in view of estimation purposes.
 While we are aware that estimation-based priors should not be routinely used in model determination,
  we remark that our reference priors  are conjugate (and thus decompose into local blocks
  precisely like the likelihood) and that they are proper. These attractive features make them natural candidates
  also
for model comparison, e.g. {\em via} Bayes factors, at least for a
preliminary and informal evaluation.

\begin{center}

\textbf{Acknowledgment}
\end{center}
The research of G.C. was partially supported by MIUR of Italy (PRIN 2005132307) and the
University of Pavia. H. M. thanks NSERC whose generous support  has made this work possible.

\section{References}

\begin{description}

\item AGRESTI, A. (2002). {\em Categorical Data Analysis}. 2nd ed. Chichester: Wiley

\vspace{-1ex}

\item ASMUSSEN, S. \& EDWARDS, D. (1983). Collapsibility and
response variables
    in contingency tables. \emph{Biometrika} \textbf{70}, 567-78.

\vspace{-1ex}

\item
BARNDORFF-NIELSEN, O.E. (1978),
{\em Information and Exponential Families in Statistical Theory}.
Chichester: Wiley.

\vspace{-1ex}

\item
BARNDORFF-NIELSEN, O.E., \& KOUDOU, A.E. (1995).
Cuts in natural exponential families. {\em Th.  Prob.
 Appl.} {\bf 40}, 361-372.

\vspace{-1ex}

\item BERGER, J. O. (2000).
Bayesian analysis: a look at today and thoughts of tomorrow.
{\em J. Am. Statist. Assoc.} {\bf 95}, 1269-1276.

\vspace{-1ex}

%

%

\item
BERGSMA, W. \& RUDAS, T. (2002). Marginal models for categorical data.
    {\em Ann. Statist.}, {\bf 30}, 140-159.

\vspace{-1ex}

%
%

\item BERNARDO, J. M. (2005). Reference analysis. In {\em Handbook of Statistics},
\textbf{25}. Eds D. K. Dey \& C. R. Rao, pp. 17-90. Amsterdam: Elsevier.

\vspace{-1ex}

%
%

%

\item
CASELLA, G. (1996).
Statistical inference and Monte Carlo algorithms (with discussion).
\emph{Test}. {\bf 5}, 249-279.

\vspace{-1ex}

%
%

%
%

\item CONSONNI, G., VERONESE, P., \& GUTIERREZ-PENA E. (2004).
Reference priors for natural exponential families having a simple quadratic variance function.
\emph{J. Multivariate Analysis}, \textbf{88}, 335-364.

\vspace{-1ex}

\item CONSONNI, G. \& LEUCARI, V. (2006).
Reference priors for discrete graphical models. \emph{Biometrika}, \textbf{93}, 23-40.

\vspace{-1ex}

\item COWELL, R.G., DAWID, A.P., LAURITZEN, S.L., \&
SPIEGELHALTER, D.J. (1999). {\em Probabilistic Networks and Expert
Systems}, New York: Springer Verlag.

\vspace{-1ex}

%
%



\item DAWID, A. P. \& LAURITZEN, S. L. (1993). Hyper Markov laws in
    the statistical analysis of decomposable models. \emph{Ann.
    Statist.} \textbf{21}, 1272-317.
%

\vspace{-1ex}

\item DELLAPORTAS, P. \& FORSTER, J.J. (1999). Markov chain Monte Carlo model determination
for hierarchical and graphical log-linear models. \emph{Biometrika}, \textbf{86}, 615-633.

\vspace{-1ex}

\item
DELLAPORTAS, P. \& TARANTOLA, C. (2005).
    Model determination for categorical data with factor level merging.
    \emph{J. R.  Statist. Soc. B}, \textbf{67}, 269-283.

\vspace{-1ex}

\item
GLONEK, G. F. V. \& MCCULLAGH, P. (1995). Multivariate logistic models. \emph{J. R. Statist. Soc. B},
 \textbf{57}, 533-546.

\vspace{-1ex}


%
%

\item
KING, R. \& BROOKS, S.P. (2001). Prior induction for log-linear models for general contingency table analysis.
\emph{Ann. Statist.}, \textbf{29}, 715-747.

\vspace{-1ex}

\item KOTZ, S., BALAKRISHNANA, N. \& JOHNSON, N.L. (2000). \textit{Continuous Multivariate
Distributions. Models and Applications}. {\bf 1}. 2nd ed. Wiley: New York.

\vspace{-1ex}

\item  LANG, J. B. and AGRESTI, A. (1994). Simultaneously modelling joint and marginal distributions of
multivariate categorical responses. \emph{J. Amer. Statist. Assoc.}, \textbf{89}, 625–632.

\vspace{-1ex}

 \item LAURITZEN, S. L. (1996). \emph{Graphical Models}. Oxford: Oxford University Press.

\vspace{-1ex}
\item LIU, J. \& MASSAM. H. (2007). The conjugate prior for discrete hierarchical log-linear models.
Submitted.

\vspace{-1ex}

\item
MADIGAN, D. \& YORK, J. (1995).
Bayesian graphical models for discrete data.
\emph{Int. Statist. Rev}, \textbf{63}, 215-232.

\vspace{-1ex}




\end{description}

\appendix
\section{Appendix}
\noindent {\bf Proof of (\ref{downwards}).}
Consider $D\subseteq C_l, D\cap R_l\not =\emptyset$ for some $l\in \{1,\ldots,k-1\}$ such that also $D\subseteq S_j$ for some $j\in \{>l\}$, then
\begin{eqnarray}
p^{C_l}(i_D)&=&\sum_{L\subseteq_0 C_l^c}p(i_D,j_L)=\sum_{L\subseteq_0 C_l^c} \exp \big(\sum_{E\subseteq_0D}\theta(i_E)+\sum_{E\subseteq_0D,G\subseteq L, j_G\in {\cal I}^*_G}\theta(i_E,j_G)\big)\nonumber\\
&=& \Big(\exp \sum_{E\subseteq_0D}\theta(i_E)\Big)\Big(1+\sum_{L\subseteq C_l^c}\exp \big(\sum_{E\subseteq_0D,G\subseteq L, j_G\in {\cal I}^*_G}\theta(i_E,j_G)
\Big)\nonumber\\
\log p^{C_l}(i_D)&=&\sum_{E\subseteq_0D}\theta(i_E)+\log \Big(1+\sum_{L\subseteq C_l^c}\exp \big(\sum_{E\subseteq_0D,G\subseteq L, j_G\in {\cal I}^*_G}\theta(i_E,j_G)
\Big)\nonumber\\
\sum_{E\subseteq_0D}\theta(i_E)&=&\log p^{C_l}(i_D)-\log \Big(1+\sum_{L\subseteq C_l^c}\exp \big(\sum_{E\subseteq_0D,G\subseteq L, j_G\in {\cal I}^*_G}\theta(i_E,j_G)
\Big)\nonumber
\end{eqnarray}
This last equality is of the form
\begin{equation}
\sum_{E\subseteq_0 D}\psi(E)=\phi(D)
\end{equation}
and therefore by Moebius inversion formula, we have
\begin{equation}
\label{72}
\psi (D)=\sum_{F\subseteq_0 D}(-1)^{|D\setminus F|}\phi(F)\;.
\end{equation}
For $l=2,\ldots,k$, let $C_{<l}=H_{l-1}\setminus C_l$. Then (\ref{72}) can be written as
\begin{eqnarray}
\theta (i_D)&=&\sum_{F\subseteq_0 D}(-1)^{|D\setminus F|}\log p^{C_l}(i_F)-\sum_{F\subseteq_0 D}(-1)^{|D\setminus F|}
\log \Big(1+\sum_{L\subseteq C_l^c}\exp \big(\sum_{H\subseteq_0F,G\subseteq L, j_G\in {\cal I}^*_G}\theta(i_H,j_G)\big)\Big)\nonumber\\
&=&\theta^{C_l} (i_D)-\sum_{F\subseteq_0 D}(-1)^{|D\setminus F|}
\log \Big(1+\sum_{L\subseteq C_l^c}\exp \big(\sum_{H\subseteq_0F,G\subseteq L, j_G\in {\cal I}^*_G}\theta(i_H,j_G)\big)\Big)\nonumber\\
&=&\theta^{C_l} (i_D)\nonumber\\
&&-\sum_{F\subseteq_0 D}(-1)^{|D\setminus F|}
\log \Big(1+\left(\sum_{L\subseteq C_{<l}}+\sum_{L\subseteq C_{>l}}+\sum _{L\subseteq C_{<l}\cup C_{>l}}\right)
\exp \big(\sum_{H\subseteq_0F,\atop{G\subseteq L,\atop  j_G\in {\cal I}^*_G}}\theta(i_H,j_G)\big)\Big)\nonumber\\
&=&\theta^{C_l} (i_D)\nonumber\\
&&-\sum_{F\subseteq_0 D}(-1)^{|D\setminus F|}
\log \Big(1+\sum_{L\subseteq C_{<l}}\exp \big(\sum_{H\subseteq_0F,\atop{G\subseteq L,\atop  j_G\in {\cal I}^*_G}}\theta(i_H,j_G)\big)\Big)
\Big(1+\sum_{L\subseteq C_{>l}}\exp \big(\sum_{H\subseteq_0F,\atop{G\subseteq L,\atop  j_G\in {\cal I}^*_G}}\theta(i_H,j_G)\big)\Big)\nonumber\\
&=&\theta^{C_l} (i_D)-\sum_{F\subseteq_0 D}(-1)^{|D\setminus F|}\log \Big(1+\sum_{L\subseteq C_{<l}}\exp \big(\sum_{H\subseteq_0F,\atop{G\subseteq L,\atop  j_G\in {\cal I}^*_G}}\theta(i_H,j_G)\big)\Big)\nonumber\\
&&\hspace{3cm}-\sum_{F\subseteq_0 D}(-1)^{|D\setminus F|}\log \Big(1+\sum_{L\subseteq C_{>l}}\exp \big(\sum_{H\subseteq_0F,\atop{G\subseteq L,\atop  j_G\in {\cal I}^*_G}}\theta(i_H,j_G)\big)\Big)\nonumber
\end{eqnarray}
We now want to show that the term
\begin{equation}
\label{t=0}
\sum_{F\subseteq_0 D}(-1)^{|D\setminus F|}\log \Big(1+\sum_{L\subseteq C_{<l}}\exp \big(\sum_{H\subseteq_0F,\atop{G\subseteq L,\atop  j_G\in {\cal I}^*_G}}\theta(i_H,j_G)\big)\Big)
\end{equation}
in the equation above is equal to zero.

Let $F$ be an arbitrary subset of $D$ and let $I=F\cap H_{l-1}$. Since $G\subseteq L \subseteq C_{<l}$, in order for $\theta(i_H,j_G),\;H\subseteq_0 F, j_G\in {\cal I}^*_G$ to be non zero, it is necessary that $H\subseteq_0 I$ and therefore
\begin{equation}
\label{equali}
\Big(1+\sum_{L\subseteq C_{<l}}\exp \big(\sum_{H\subseteq_0F,\atop{G\subseteq L,\atop  j_G\in {\cal I}^*_G}}\theta(i_H,j_G)\big)\Big)=\Big(1+\sum_{L\subseteq C_{<l}}\exp \big(\sum_{H\subseteq_0I,\atop{G\subseteq L,\atop  j_G\in {\cal I}^*_G}}\theta(i_H,j_G)\big)\Big)\;.
\end{equation}
We see that the right hand side of (\ref{equali}) above is the same for all $F\subseteq_0 D$ that have the same intersection $I$ with $H_{l-1}$.
We therefore consider all such $F$'s.  Since $D\cap R_l\not=\emptyset$, there are as many such $F$'s with $|D\setminus F|$ odd as there are with  $|D\setminus F|$ even and therefore  from (\ref{equali}), it follows that, for a given $I$,
\begin{equation}
\label{equal0}
\sum_{F\subseteq_0 D\atop F\cap H_{l-1}=I}(-1)^{|D\setminus F|}\log \Big(1+\sum_{L\subseteq C_{<l}}\exp \big(\sum_{H\subseteq_0F,\atop{G\subseteq L,\atop  j_G\in {\cal I}^*_G}}\theta(i_H,j_G)\big)\Big)=0\;.
\end{equation}
Since this is true for all $I\subseteq_0D\cap H_{l-1}$, it follows immediately from (\ref{equal0}) that (\ref{t=0}) is equal to zero and we have
\begin{eqnarray}
\theta (i_D)&=&\theta^{C_l} (i_D)-\sum_{F\subseteq D}(-1)^{|D\setminus F|}\log \Big(1+\sum_{L\subseteq C_{\{>l\}}}\exp \big(\sum_{H\subseteq_0F,\atop{G\subseteq L,\atop  j_G\in {\cal I}^*_G}}\theta(i_H,j_G)\big)\Big)\nonumber
\end{eqnarray}
Since $D\subseteq S_j\cap C_l$ for some $m\in \{>l\}$ and
$G\subseteq L\subseteq C_m\setminus C_l$ is non empty, in the right
hand side of the equation above, we have that either
$\theta(i_H,j_G)=\theta^{C_m}(i_H,j_G)$ or that $\theta(i_H,j_G)$
can be expressed using (\ref{downwards}) recursively and therefore
$\theta(i_D)$ can be expressed in terms of $\theta^{C_m}(i_E), m\in
\{>l\}, E\subseteq C_m, i_E\in {\cal I}^*_E$. Formula
(\ref{downwards}) is thus proved.
\end{document}